\documentclass[aps,prl,floatfix, amsmath,preprint, showpacs,showkeys,superscriptaddress]{revtex4-1}
\usepackage{graphicx,xcolor}
\usepackage{amssymb,amsthm}
\usepackage{epstopdf}
\usepackage{amsfonts}
\usepackage{amsmath}
\usepackage{xcolor}
\usepackage{amssymb}
\usepackage{epstopdf}
\usepackage{dsfont}

\newcommand{\bE}{\mathbb{E}}

\begin{document}

\title{Fast and slow domino regimes in transient network dynamics}

\author{Peter Ashwin}
\affiliation{Department of Mathematics, University of Exeter, Exeter EX4 4QF, UK and \\EPSRC Centre for Predictive Modelling in Healthcare, University of Exeter, Exeter, EX4 4QJ, UK.}
\author{Jennifer Creaser}
\affiliation{Department of Mathematics, University of Exeter, Exeter EX4 4QF, UK and \\EPSRC Centre for Predictive Modelling in Healthcare, University of Exeter, Exeter, EX4 4QJ, UK.}
\author{Krasimira Tsaneva-Atanasova}
\affiliation{Department of Mathematics and Living Systems Institute, University of Exeter, Exeter EX4 4QF, UK and \\EPSRC Centre for Predictive Modelling in Healthcare, University of Exeter, Exeter, EX4 4QJ, UK.}

\begin{abstract}
It is well known that the addition of noise to a multistable dynamical system can induce random transitions from one stable state to another. For low noise, the times between transitions have an exponential tail and Kramers' formula gives an expression for the mean escape time in the asymptotic limit. If a number of multistable systems are coupled into a network structure, a transition at one site may change the transition properties at other sites. We study the case of escape from a ``quiescent'' attractor to an ``active'' attractor in which transitions back can be ignored. There are qualitatively different regimes of transition, depending on coupling strength. For small coupling strengths the transition rates are simply modified but the transitions remain stochastic. For large coupling strengths transitions happen approximately in synchrony - we call this a ``fast domino" regime. There is also an intermediate coupling regime some transitions happen inexorably but with a delay that may be arbitrarily long - we call this a ``slow domino" regime. We characterise these regimes in the low noise limit in terms of bifurcations of the potential landscape of a coupled system. We demonstrate the effect of the coupling on the distribution of timings and (in general) the sequences of escapes of the system.
\end{abstract}



\maketitle

A number of important physical, biological and socio-economic questions involve understanding how a dynamical change of one subsystem within a network affects other subsystems that are coupled to it. Indeed, there is extensive work on noisy coupled bistable units, motivated by trying to understand collective response and phase transitions. This includes work on stochastic resonance on networks \cite{intro1,intro2}. For example, \cite{intro11} uses a master equation approach while \cite{intro12,intro13} consider noise-induced switching of bistable nodes in complex networks. Much of this work aims to explain properties of attracting (statistically steady) states perturbed by noise; nonetheless, many important questions are related to the transient dynamics of networks affected by noise.

We consider transient noise-induced behaviour in a network of asymmetric bistable attractor systems, where noise induces an effectively irreversible transition spread through coupling. Each node (corresponding to a subsystem) is assumed to have two states, a  shallow marginally stable mode (the ``quiescent'' state) and a deep more stable mode (the ``active'' state) that is consequently more resistant to noise. We start with the system in the marginally stable mode  and say it ``escapes" when it crosses some threshold to the deeply stable mode. The time of first escape is a random variable that is jointly determined by the nonlinear dynamics and the noise process. The assumption of asymmetry means that escape from the deeper state occurs very rarely and so we can view the process as an irreversible cascade of escapes, similar to a cascade of toppling dominos. The coupling of the systems can promote (or hinder) escape of others on the network and may cause certain sequences of escape to appear preferentially depending on coupling strength. In this paper we highlight that the timings and sequences of escapes are effectively ``emergent properties'' of the system, and we demonstrate that these properties can be usefully classed by coupling strength into qualitatively different regimes.

We consider an idealization of behaviour that has been seen in a variety of applications: this includes (a) signal propagation by sequential switching between asymmetric stable states (observed experimentally in chains of bistable electronic circuits \cite{intro3} or in cases where the bistability is noise-induced \cite{intro14}) (b) waves along unidirectionally coupled chains (or lattices) of bistable nodes with forcing at one end \cite{intro4} (c) photoinduced phase transitions in spin-crossover materials with bistable dynamic potentials \cite{intro5,intro6,intro7} (d) avalanches of gene activation in gene regulatory pathways to drive cell differentiation/development/cancer \cite{intro8,intro9} (e) cell fate in biofilm formation \cite{intro10}. Other applications that could benefit from a better understanding of similar transient dynamics induced by noise include (a) the contagion of bank defaults in a system of financial institutions interconnected by mutual loans \cite{GaiKap10,HalMay11,chinazzi13,summer13}, (b) interconnections between ``tipping elements'' \cite{AshWieVitCox2012,Lenton_etal_2008}, (c) the role of spreading of abnormal large-amplitude oscillators in modelling onset of epileptic seizures \cite{kalitzin10,benj12pheno} (d) multiple organ failure \cite{Parkeretal2010} or (e) cascading failures in power systems \cite{Dobsonetal2007}.

The role of coupling strength in noise-induced transitions on networks is considered by \cite{BFG2007a,BFG2007b} for idealised symmetric bistable systems. Neiman \cite{Neiman_1994} shows similar synchronization effects in coupled stochastic bistable systems and \cite{Mateos_Alatriste_2010} in coupled ratchet systems. The authors of \cite{BFG2007a,BFG2007b} give rigorous mathematical results that identify the existence of different regimes of synchronization of escapes in the low noise limit that can be linked to changes in the structure of underlying system attractors (see for example \cite{ChaFer05} for some review of the role of coupling in the noise-free context).  In particular, \cite{BFG2007a} identify that the most likely sequences of escape and how their probabilities change qualitatively with coupling strength: there can be synchronized transitions in the strong coupling limit. Many properties  of the transitions can be understood using Friedlin-Wentzell methodology and the Eyring-Kramers formula \cite{BG2006,Berglund2013} to study the pathwise properties of transitions between attractors.

We show in the context of asymmetric potentials that there are typically several qualitatively different regimes in the transient sequences of escapes. These regimes of weak, intermediate and strong coupling, and the intermediate case may be quite complicated, but in general there are qualitative changes in behaviour for the weak noise limit that can be characterised in terms of bifurcations of steady states of the noise-free system. As a row of toppling dominos depends on the properties and spacing of the dominos~\cite{lee10domino},  we identify different domino effects that can be characterised by different coupling regimes. Specifically, we identify ``slow domino" and ``fast domino" regimes corresponding to intermediate and strong coupling regimes, respectively. Within these different regimes, certain sequences of escape may be preferred by the coupling, and the distribution of times to next escape may have significant deviations from exponential.

We consider a diffusively coupled network of prototypical asymmetric bistable nodes under the influence of additive noise for an asymmetric case of the Schl\"{o}gl model \cite{Malchow_etal_1983}. For $N=2$ nodes and bidirectional coupling there are qualitative changes in the escape time distributions as the coupling strength increases~\cite{frank82stoch}. For $N=3$ nodes with unidirectional coupling, we show that, although the mean and distributions of escape times of an individual node are not much affected by the coupling, the probability of a given sequence appearing and the distribution of timings within the sequence of escapes can be greatly affected. 

We consider a network where each node is governed by a bistable system
\begin{equation}
\label{eq:onenode}
\dot{x}=f(x,\nu):=-(x-1)(x^2-\nu)
\end{equation}
so that $f=-V'(x)$ with potential $V(x)=\frac{1}{4}x^4-\frac{1}{3}x^3+\nu(x-\frac{1}{2}x^2)$. We suppose that nodes are coupled into a network and subjected to additive noise. For $0<\nu\ll 1$ the stable states are not interchangeable by any symmetry: there is a quiescent attractor at $x=x_Q:=-\sqrt{\nu}$ and an active attractor  at $x=x_A:=1$; there is an unstable separating equilibrium at $x=x_S:=\sqrt{\nu}$. Stationary distributions of this model are examined in \cite{Malchow_etal_1983}. For nodes $i=1,\cdots,N$ the network is assumed to evolve according to the SDE
\begin{equation}
dx_i=\bigl[f(x_i,\nu)+\beta\sum_{j\in N_i} (x_j-x_i)\bigr]dt + \alpha \,dw_i
\label{eq:Nnode}
\end{equation}
where $N_i$ are the neighbours that provide inputs to node $i$, $\beta$ is the coupling strength, $\alpha$ the strength of the additive noise and $w_i$ are independent Wiener processes. 

\begin{figure}	\includegraphics[width=0.6\textwidth]{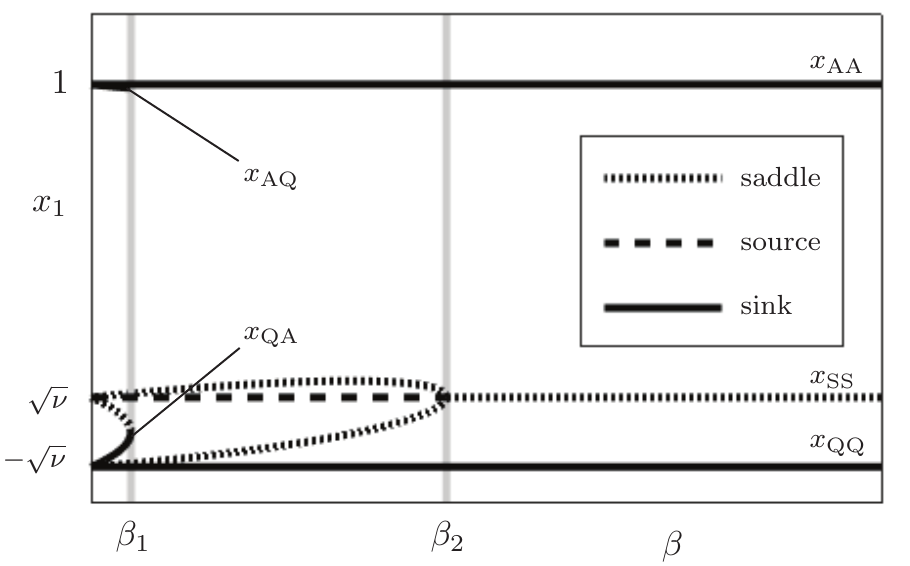}
\caption{
	Bifurcation diagram for the system of two bidirectionally coupled nodes (\ref{eq:bistable2node}) with $\alpha=0$ and $\nu=0.01$ projected into the $(\beta,x_1)$ plane, where $\beta$ is the coupling strength (cf \cite[Fig 2]{frank82stoch}). We are interested in how the system escapes from the quiescent attracting state $x_{QQ}$ to the active attracting state $x_{AA}$ under the influence of low-amplitude noise, $0<\alpha\ll 1$. The three regimes that exist in terms of the structures that must be overcome for the transition have parallels in more general cases. In this case they are divided by a saddle-node (fold) bifurcation at $\beta_1=0.0101$ and a pitchfork bifurcation of the separating saddles at $\beta_2=0.09$. In the weak coupling regime $\beta<\beta_1$ the escape will be via an additional attractor, $x_{QA}$ or $x_{AQ}$, while in the strong coupling (``fast domino'' regime) $\beta>\beta_2$, the escapes are approximately synchronised and pass near $x_{SS}$. Escapes in the intermediate coupling (``slow domino'' regime) $\beta_1<\beta<\beta_2$ are associated with escape over a symmetry broken saddle.
\label{fig:bistable_2nodebif}
}
\end{figure}

In the case $N=2$ with bidirectional coupling \cite{frank82stoch} we have
\begin{align}
dx_1&=\left[f(x_1,\nu)+\beta(x_2-x_1)\right]dt + \alpha\, dw_1,\label{eq:bistable2node}\\
dx_2&=\left[f(x_2,\nu)+\beta(x_1-x_2)\right]dt + \alpha\, dw_2\nonumber
\end{align}
where in the noise-free case $\alpha=0$ there are equilibria at $x_{QQ}:=(x_Q,x_Q)$, $x_{SS}:=(x_S,x_S)$ and $x_{AA}:=(x_A,x_A)$ for any $\beta$. Up to six more equilibria depend on $0\leq \beta$ and $0<\nu<1$. The regimes noted in \cite{frank82stoch} can be precisely characterized: one can verify that the number of solutions changes at a saddle node bifurcation when
$$
-27\beta^3+(27\nu+9)\beta^2-9(\nu+{\textstyle \frac{1}{3}})^2\beta+\nu(\nu-1)=0.
$$
For small $\nu$ this implies there is a saddle-node for $\beta=\beta_1>0$. A pitchfork bifurcation occurs at intermediate $\beta_2=(\sqrt{\nu}-4\nu+3\nu^{3/2})/(1-3\sqrt{\nu})$. Let $x_{QS}$ denote the branch of equilibria that continues from $(x_Q,x_S)$ at $\beta=0$. We note $x_{SA}$ (saddle) and $x_{QA}$ (stable) meet while simultaneously $x_{AS}$ (saddle) and $x_{AQ}$ (stable) meet at the saddle-node at $\beta_1$. The branches $x_{QS}$ and $x_{SQ}$ meet $x_{SS}$ at the pitchfork bifurcation at $\beta_2$. Observe that there are three qualitatively different regimes of coupling depending on whether there are nine ($\beta<\beta_1$), five ($\beta_1<\beta<\beta_2$) or three ($\beta>\beta_2$) equilibria.  The bifurcation diagram for $\nu=0.01$ is shown in Figure~\ref{fig:bistable_2nodebif}:  in this case $\beta_1=0.0101$ and $\beta_2=0.09$. 

%%%%%%%%%%%%%%%%%%%%%%%%%%%%%%%%%%
%%%%%%%%%%%%%%%%%%%%%%%%%%%%%%%%%%

We give initial condition $x_i(0)=x_Q$ for (\ref{eq:Nnode}) and pick a threshold $x_S<\xi<x_A$. The first escape time of node $i$ is the random variable $\tau^{(i)}=\inf\{t>0~:~x_i(t)>\xi\}$ that depends on the network, the parameters and the particular noise path: it has a distribution implied by that of the noise. Independence of the $w_i$ means that (with probability one) no two escapes will occur at the same time and so we can assume there is a permutation $s(i)$ of $\{1,\ldots,N\}$ such that $\tau^{s(i)}<\tau^{s(j)}$ for any $i<j$. We denote by $\mathbb{P}(s)$ the probability of a sequence $s$ being realised and define the time of the $i$th escape by $\tau^{i}=\tau^{s(i)}$: we use the convention $\tau^{0}=0$. The time between escapes $j$ and $k>j$ is denoted $\tau^{k|j}=\tau^{k}-\tau^{j}$, with means $T^{(i)}=\bE[\tau^{(i)}]$ and $T^{k|j}=\bE[\tau^{k|j}]$. Note that for $\beta=0$ all sequences are equally likely, meaning $\mathbb{P}(s)=1/N!$. 

In networks of the form (\ref{eq:bistable2node}), as long as $0<\nu<1$ so that $x_Q$ is linearly stable, the $\tau^{(i)}$ are independent random variables with exponential tails for $\beta=0$ whose mean can be approximated using the one-dimensional Kramers' formula (e.g. \cite{Berglund2013}) which states in the limit $\alpha\rightarrow 0$:
\begin{equation}
\label{eq:Kramersb0}
T^{(i)} \approx \frac{2\pi}{\sqrt{V''(x_Q)|V''(x_S)|}} \text{e}^{ \frac{2}{\alpha^2}[V(x_S)-V(x_Q)]}.
\end{equation}
We show that the distributions $\tau$ and $\mathbb{P}(s)$ change in subtle ways on increasing $\beta$. 

Persistence of the hyperbolic fixed points and robustness of connections means there is a weak coupling regime: for small enough $\beta>0$, the quiescent states are perturbed but not destroyed, and escape of one node modifies the rate of escape of the other nodes. However the means (\ref{eq:Kramersb0}) should vary continuously with the parameter.
For the strong coupling (synchronized) regime \cite{Neiman_1994,BFG2007a}: for large $\beta$ the nodes synchronize and there is strong dependence, meaning they escape {\em en masse}: hence ``fast domino". For the intermediate coupling regime where escape of one node leads to a delayed (but essentially deterministic) response from the other units: hence ``slow domino".

\begin{figure}	
\includegraphics[width=0.6\textwidth]{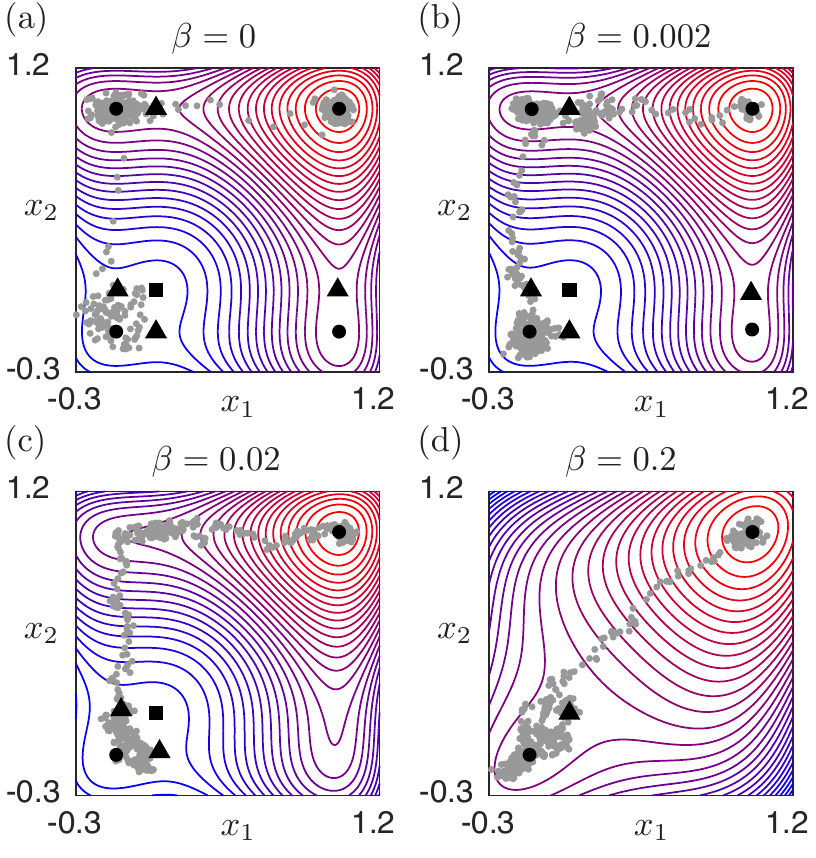}
\caption{
Level sets of $\tilde{V}$ for $N=2$ bidirectionally coupled nodes (\ref{eq:bistable2node}) with fixed $\nu=0.05$ and four values of $\beta$. The equilibria for $\alpha=0$ are marked as $\bullet$   sinks, $\blacksquare$   sources and $\blacktriangle$  saddles. Typical noise paths starting at $x_{QQ}$ are shown in each panel computed for (\ref{eq:bistable2node}) and for $\alpha=0.1$. The panels show typical escapes of (a) uncoupled (b) weakly coupled (c) intermediate coupled (``slow domino'') and (d) strongly coupled (``fast domino'') regimes.
\label{fig:bistable_2node_noise}
}
\end{figure}

We illustrate these differences for (\ref{eq:bistable2node}) in Figure~\ref{fig:bistable_2node_noise}, which shows the behaviour of escapes from $x_{QQ}$ in the weak noise limit with $\nu=0.01$ fixed and depending on $\beta$, where the SDE is solved using a fixed timestep Heun method.  The symmetry in the coupling of the system can be seen as a reflection about the line $x_1=x_2$.
The coupled system (\ref{eq:bistable2node}) can be seen as a noise perturbed potential flow for
$\tilde{V}(x_1,x_2)= V(x_1)+V(x_2)+\frac{1}{2}\beta(x_1-x_2)^2$ (we suppress the $\nu$ and $\beta$ dependence).
The mean escape time between two minima of the potential can be estimated using a multidimensional Kramers' formula: the mean time from $x^*$ to $y^*$ over the minimum height pass saddle (`gate') at $z^*$ is
$$
T(x^*,z^*,y^*) \approx P(x^*,z^*) e^{\frac{2}{\alpha^2} [\tilde{V}(z^*)-\tilde{V}(x^*)]}
$$
for $\alpha\rightarrow 0$, where the prefactor $P$ depends on the Hessian $\nabla^2\tilde{V}(z^*)$ (see e.g. \cite{Berglund2013}).  Note that to this leading order $T$ is independent of $y^*$. 

%%%%%%%%%%%%%%%%%%%%%%%%%%%%%%%%%%
%%%%%%%%%%%%%%%%%%%%%%%%%%%%%%%%%%

We estimate the dependence of mean time $T^{2|0}=T^{2|1}+T^{1|0}$ of escape for (\ref{eq:bistable2node}) on coupling, where there may be multiple paths of escape. If $\widetilde{T}(x^*,\tilde{z}^*, y^*)$ is the mean time of escape assuming it takes path $\tilde{z}^*$ out of $G$ possible symmetrically equivalent gates, then $\widetilde{T}(x^*,\tilde{z}^*,y^*)=\frac{1}{G}T(x^*,z^*,y^*)$, where $z^*$ is associated with multiple paths of escape.

In the {\em weak coupling regime} $0<\beta<\beta_1$ each symmetric path is equally probable and so $2T^{1|0} \approx  \widetilde{T}(x_{QQ},x_{QS},x_{QA})+\widetilde{T}(x_{QQ},x_{SQ},x_{AQ})$,
 while $2T^{2|1}\approx T(x_{QA},x_{SA},x_{AA})+T(x_{AQ},x_{AS},x_{AA})$. Hence
\begin{align}
T^{2|0}\approx  {\textstyle \frac{1}{2}}T(x_{QQ},x_{QS},x_{QA})+T(x_{QA},x_{SA},x_{AA}).
\end{align}

In the {\em intermediate coupling regime} (``slow domino" regime) $\beta_1<\beta<\beta_2$ there is a one-step escape process, but there are two possible gates that can be traversed:
\begin{align}
T^{2|0}\approx {\textstyle \frac{1}{2}}\left[T(x_{QQ},x_{SQ},x_{AA})+T(x_{QQ},x_{QS},x_{AA})\right].
\end{align}
Note that this asymptotic expression will be non-uniform in $\beta$: near $\beta=\beta_1$ there will be a long deterministic delay associated with passage past the region of the saddle-node as is evident in Figure~\ref{fig:bistable_2node_noise}(c).

In the {\em strong coupling regime} (``fast domino" regime) $\beta>\beta_2$ there is a one-step escape process with a unique gate:
\begin{align}
T^{2|0}\approx T(x_{QQ},x_{SS},x_{AA}).
\end{align}
Each of these regimes will give a different scaling in the limit $\alpha\rightarrow 0$, while the scalings at crossovers between regimes are accessible to generalizations of Kramers' formula for passage over nonhyperbolic saddles \cite{Berglund2013}. This is explored in more detail in \cite{CTA}, including  computing the timing of the escape once the gate has been traversed in the intermediate and strong coupling regimes.

\begin{figure*}[ht]
\includegraphics[width=16.5cm]{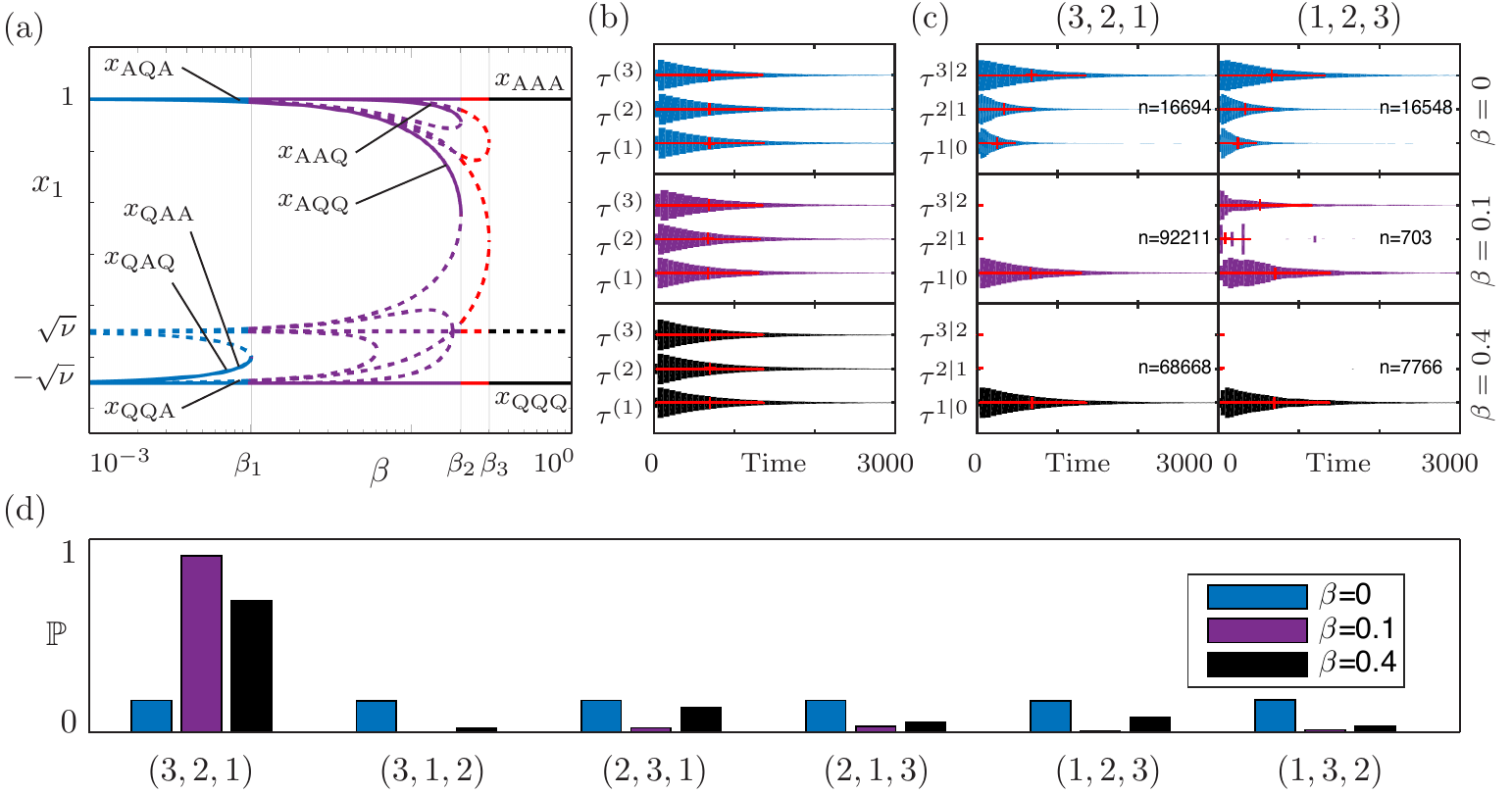}
\caption{
(a) Bifurcation diagram showing $x_1$ vs $\beta$ (log axis) for \eqref{eq:bistable3chain} with $\nu=0.01$ and no noise $\alpha=0$: dashed branches are unstable. In the weak coupling regime ($\beta<\beta_1=0.0101$, blue) all branches continue from $\beta=0$. There are two intermediate (slow domino) coupling regimes: for the lower one ($\beta_1<\beta<\beta_2\approx0.2025$, purple) there are still stable and unstable partially escaped states while for ($\beta_2<\beta<\beta_3\approx0.3035$, red) there are only partially escaped saddles. For the strong (fast domino) coupling regime $\beta>\beta_3$ all equilibria are synchronized in the absence of noise. For (b)--(d)  we computed $10^5$ samples using $\alpha=0.03$ for $\beta=0$ (blue), $0.1$ (purple) and $0.4$ (black). Panel (b) shows violin plots of the distribution of escape times $\tau^{(i)}$ of node $i$: observe that these change little with coupling. The red cross indicates mean (vertical) and $+/-$ one standard deviation (horizontal).  Panel (c) shows the distribution of sequential escape times $\tau^{k|k-1}$ for $k=1,2,3$, for sequences $(3,2,1)$ and $(1,2,3)$. The number of samples $n$ (out of $10^5$) that undergo this sequence of escapes is shown. Panel (d) shows the probability of a given sequence being realised. In the strongly coupled case $\beta=0.4$ the escapes are almost always synchronized, and the most frequent sequence is $(3, 2, 1)$. The case $\beta=0.1$ and sequence $(1,2,3)$ is an example of a non-synchronous escape in the intermediate coupling regime; the third escape typically occurs some time after the first two: see Table~\ref{tab:vals}. 
\label{fig:bistable_3chain}
}	
\end{figure*}

For a more general network, the sequence of escapes of the network depends not only on the number of nodes that have already escaped but also the sequence in which they escape. We consider a unidirectionally coupled chain of $N=3$ bistable systems (\ref{eq:Nnode}) where the input sets $N_i$ for node $i$ are given by $(N_1,N_2,N_3)=(\{2\},\{3\},\{\})$:
\begin{align}
dx_1&=\left[f(x_1,\nu)+\beta(x_2-x_1)\right]dt + \alpha\, dw_1\nonumber\\
dx_2&=\left[f(x_2,\nu)+\beta(x_3-x_2)\right]dt + \alpha\, dw_2\label{eq:bistable3chain}\\
dx_3&=\left[f(x_3,\nu)\right]dt + \alpha\, dw_3.\nonumber
\end{align}
Figure~\ref{fig:bistable_3chain} illustrates the three coupling regimes; the weak coupling regime ($\beta<\beta_1$), intermediate coupling (slow domino) ($\beta_1<\beta<\beta_3$), and strong coupling (fast domino) ($\beta>\beta_3$) regimes for this system. Note that intermediate coupling can be split further into two sub-regimes at $\beta_2$. There are qualitative changes in the asymptotic behaviour of sequential escapes on changing $\beta$, with strongly synchronized escapes for strong coupling. 

%%%%%%%%%%%%%%%%%%%%%%%%%%%%%%%%%%
%%%%%%%%%%%%%%%%%%%%%%%%%%%%%%%%%%

To characterise the distribution of times of $n$th escape we consider the coefficient of variation of $\tau$ given by
$$
CV(\tau)=\sigma(\tau)/\bE[\tau]
$$
where $\sigma(\tau)$ denotes the standard deviation For $\beta=0.0$ (and for all first escapes) we have $CV(\tau^{k|k-1})\approx 1$, indicating an exponential distribution. In the intermediate coupling (slow domino) regime $\beta=0.1$ the most likely sequence is $(3,2,1)$: considering only this sequence for the data in Figure~\ref{fig:bistable_3chain} we find $CV(\tau^{1|0})=0.9608$, $CV(\tau^{2|1})=0.3308$ and $CV(\tau^{3|2})=0.2210$ - after the first (approximately exponentially distributed) escape the remaining escapes are close to deterministic ($\bE[\tau^{2|1}]=4.087$, $\bE[\tau^{3|2}]=4.797$). On the other hand, for a rarer sequence $(1,2,3)$ in the intermediate regime we find $CV(\tau^{1|0})=0.9783$, $CV(\tau^{2|1})=3.662$ and $CV(\tau^{3|2})=1.27$ - after the first exponentially distributed escape there are very large variations in escape time. Finally, in the strongly coupling (fast domino) regime $\beta=0.4$ and the most likely sequence $(3,2,1)$ we have  $\bE[\tau^{2|1}]=0.6568$, $\bE[\tau^{3|2}]=0.9664$.  Table~\ref{tab:vals} gives the probability, mean and coefficient of variation for sequential escape times of the simulations shown in Figure~\ref{fig:bistable_3chain}. Note that as $\beta$ increases, the system remains closer to synchronization, leading to an increasing randomization of the sequence of escapes caused by fluctuations about the synchronized state.

%TTTTTTTTTTTTTTTTTTTTTTTTTTTTTTTT
\begin{table}[ht]
	\caption{Data Table. For the simulations shown in Figure~\ref{fig:bistable_3chain}, the columns in this table show the sequence of escape, the probability $\mathbb{P}$ that a sequence will be realised, followed by the mean, standard deviation and coefficient of variation of $\tau^{k|k-1}$ conditional on this sequence for $k=1,2,3$. }
	\label{tab:vals}
	\centering 
	\renewcommand{\arraystretch}{1.1}
	\begin{tabular}{ | c | c || c | c | c | c || c |  c | c | c || c | c | c | c |}
    \multicolumn{14}{c}{\textbf{$\beta$ = 0}: Uncoupled systems }\\ \hline 
		Sequence & $\mathbb{P}$ & $\tau$ & $\bE(\tau)$ & $\sigma(\tau)$ & $CV(\tau$) & $\tau$ & $\bE(\tau)$ & $\sigma(\tau)$ & $CV(\tau$) & $\tau$ & $\bE(\tau)$ & $\sigma(\tau)$ & $CV(\tau$)  \\ \hline \hline
		(3, 2, 1) & 0.167 & $\tau^{1|0}$ & 244.53 & 221.98 & 0.91 & $\tau^{2|1}$ & 334.87 & 340.60 & 1.02 & $\tau^{3|2}$ & 673.07 & 668.26 & 0.99 \\  \hline
		(3, 1, 2) & 0.166 & $\tau^{1|0}$ & 245.94 & 222.72 & 0.91 & $\tau^{2|1}$ & 333.61 & 330.46 & 0.99 & $\tau^{3|2}$ & 662.49 & 661.12 & 1.00 \\ \hline 
		(2, 3, 1) & 0.167 & $\tau^{1|0}$ & 246.58 & 226.22 & 0.92 & $\tau^{2|1}$ & 332.64 & 329.08 & 0.99 & $\tau^{3|2}$ & 668.02 & 674.47 & 1.01 \\ \hline
		(2, 1, 3) & 0.167 & $\tau^{1|0}$ & 243.26 & 223.67 & 0.92 & $\tau^{2|1}$ & 334.81 & 331.77 & 0.99 & $\tau^{3|2}$ & 671.92 & 665.28 & 0.99 \\ \hline 
		(1, 2, 3) & 0.165 & $\tau^{1|0}$ & 243.57 & 223.05 & 0.92 & $\tau^{2|1}$ & 337.94 & 337.15 & 1.00 & $\tau^{3|2}$ & 664.35 & 655.76 & 0.99 \\ \hline
		(1, 3, 2) & 0.168 & $\tau^{1|0}$ & 246.26 & 224.39 & 0.91 & $\tau^{2|1}$ & 329.51 & 329.09 & 1.00 &  $\tau^{3|2}$ & 667.31 & 667.83 & 1.00 \\ \hline
		 \multicolumn{14}{c}{\textbf{$\beta$ = 0.1}: Intermediate coupling regime  (``slow domino effect" ) }\\ \hline 
		(3, 2, 1) & 0.922 & $\tau^{1|0}$ & 658.98 & 633.17 & 0.96 & $\tau^{2|1}$ &  \phantom{61}4.09 &  \phantom{30}1.36 & 0.33  & $\tau^{3|2}$ & \phantom{50} 4.80  &  \phantom{67}1.06 & 0.22 \\  \hline
		(3, 1, 2) & 0.002 & $\tau^{1|0}$ & 730.13 & 658.49 & 0.90 & $\tau^{2|1}$ &  \phantom{61}2.26 &  \phantom{30}1.42 & 0.63 & $\tau^{3|2}$ &  \phantom{50}1.12 & \phantom{67}1.01 & 0.90 \\  \hline
		(2, 3, 1) & 0.024 & $\tau^{1|0}$ & 652.22 & 611.87 & 0.94 & $\tau^{2|1}$ &  \phantom{61}1.50 &  \phantom{30}1.27 & 0.85 & $\tau^{3|2}$ &  \phantom{50}2.97 & \phantom{67}1.55 & 0.52  \\  \hline
		(2, 1, 3) & 0.031 & $\tau^{1|0}$ & 666.43 & 647.67 & 0.97 & $\tau^{2|1}$ &  \phantom{61}3.54 &  \phantom{30}1.70 & 0.48& $\tau^{3|2}$ & 487.84 & 673.65 & 1.38 \\  \hline
		(1, 2, 3) & 0.007 & $\tau^{1|0}$ & 704.30 & 689.06 & 0.98 & $\tau^{2|1}$ &  \phantom{1}82.71 & 302.97 & 3.66  & $\tau^{3|2}$ & 509.47 & 647.88 & 1.27  \\  \hline
		(1, 3, 2) & 0.014 & $\tau^{1|0}$ & 703.84 & 663.34 & 0.94 & $\tau^{2|1}$ & 617.64 & 665.10 & 1.08& $\tau^{3|2}$ &  \phantom{50}3.93 & \phantom{67}1.46 & 0.37 \\  \hline 
		 \multicolumn{14}{c}{\textbf{$\beta$ = 0.4}: Strong coupling regime (``fast domino effect" ) }\\ \hline 
		(3, 2, 1) & 0.687& $\tau^{1|0}$ & 688.02 & 662.25 & 0.96 & $\tau^{2|1}$ &  \phantom{1}0.66 & \phantom{10}0.38 & 0.58   & $\tau^{3|2}$ & 0.97 & 0.40 & 0.41  \\  \hline
		(3, 1, 2) & 0.024 & $\tau^{1|0}$ & 708.41 & 691.41 & 0.98 & $\tau^{2|1}$ &  \phantom{1}0.36 &  \phantom{10}0.27 & 0.75 & $\tau^{3|2}$ & 0.21 & 0.18 & 0.86 \\  \hline
		(2, 3, 1) & 0.128 & $\tau^{1|0}$ & 690.46 & 682.03 & 0.99 & $\tau^{2|1}$ &  \phantom{1}0.29 &  \phantom{10}0.25 & 0.86  & $\tau^{3|2}$ & 0.62 & 0.39 & 0.63  \\  \hline
		(2, 1, 3) & 0.053 & $\tau^{1|0}$ & 702.68 & 681.17 & 0.97  & $\tau^{2|1}$ & \phantom{1}0.41 &  \phantom{10}0.31 & 0.76   & $\tau^{3|2}$ & 0.50 & 0.53 & 1.06 \\  \hline
		(1, 2, 3) & 0.078 & $\tau^{1|0}$ & 695.96 & 680.09 & 0.98 & $\tau^{2|1}$ &  \phantom{1}4.00 &  \phantom{1}49.62 & 12.41 & $\tau^{3|2}$ & 0.76 & 0.70 & 0.92  \\  \hline
		(1, 3, 2) & 0.030 & $\tau^{1|0}$ & 694.73 & 651.60 & 0.94 & $\tau^{2|1}$ & 17.54 & 151.01 & 8.61& $\tau^{3|2}$ & 0.30 & 0.24 & 0.80  \\  \hline 		
	\end{tabular}
\end{table}
%TTTTTTTTTTTTTTTTTTTTTTTTTTTTTTTT

For general heterogeneous networks it is still possible to classify the interactions between nodes $x_i$ and $x_j$ as weak, intermediate or strong depending on whether escape of node $x_i$ modifies the rate of noise-induced escape of $x_j$, whether $x_j$ will undergo a deterministic escape in a bounded time or whether $x_j$ will be synchronized in its escape with $x_i$, respectively. This will depend on the state of other nodes that are connected to $x_i$ and $x_j$, and so the classification of the interaction is, in general, state and sequence dependent.

The changes in distribution of timings and sequences of escapes in stochastically perturbed coupled networks can be usefully thought of as an emergent behaviour of the network. In particular, even for intermediate or strong coupling where there are no symmetry broken attractors in the noise-free case, the asymptotic behaviour of the sequence of escapes is qualitatively different in the low noise limit. A study of such sequential escapes will be of interest in a variety of situations where stochastic forcing of individual sites with asymmetric attractors interacts with the coupling strength to change the sequence of escapes.  For example, \cite{CTA} use this to explain some phenomena in the networks of coupled oscillatory bistable units considered in \cite{benj12pheno}.

\begin{acknowledgments} The authors gratefully acknowledge the financial support of the EPSRC via grant EP/N014391/1. We thank the anonymous referees for their comments, criticisms and suggestions. PA gratefully acknowledges the European Union's Horizon 2020 research and innovation programme for the ITN CRITICS under Grant Agreement number 643073 for providing opportunities to discuss this work with members of the CRITICS network.
\end{acknowledgments}

%\newpage

%%%%%%%%%%%%%%%%%%%%%%%%

\end{document}